\documentstyle[12pt]{article}
\begin{document}
\begin{center}
{\large {\bf INDUCED REPRESENTATIONS  
OF THE TWO PARAMETRIC QUANTUM DEFORMATION 
U$_{pq}$[gl(2/2)]}} 
\vskip .7truecm
{\bf Nguyen Anh Ky}\\
     
\normalsize Department of Physics, Chuo University, Kasuga, 
Bunkyo--ku\\ Tokyo 112-8551, {\bf Japan}\\[1mm]
and\\[1mm]
Institute of Physics, P.O. Box 429, Bo Ho, Hanoi 10000, 
{\bf Vietnam}
\vskip 1truecm
{\bf Abstract} \\[1mm]
\end{center}
     
 The two--parametric quantum superalgebra $U_{p,q}[gl(2/2)]$
and its induced representations are considered. A method for 
constructing all finite--dimensional irreducible representations 
of this quantum superalgebra is also described in detail.
It turns out that finite--dimensional representations of the 
two--parametric $U_{p,q}[gl(2/2)]$, even at generic deformation 
parameters, are not simply trivial deformations from those of the 
classical superalgebra $gl(2/2)$, unlike the one--parametric 
cases.\\[.5cm]

\underline{Running title} : ~ Representations of 
$U_{p,q}[gl(2/2)]$.

\vspace*{2mm}

\underline{PACS numbers} : ~ 02.20Tw, 11.30Pb.

\underline{MSC--class.} : 81R50; 17A70.
\newpage
\begin{flushleft}
{\large {\bf I. Introduction}}
\end{flushleft}
\vspace*{1mm}
     
  Introduced in 80's as a result of the study on quantum
integrable systems and Yang--Baxter equations \cite{collect} 
the quantum groups \cite{frt}--\cite{woro} have been 
intensively investigated in different aspects. Since then 
many (algebraic and geometric) structures and various 
applications of quantum (super) groups have been found
(see in this context, for example, Refs. 
\cite{2dim,majid,chari,kass}). It turns out that  
quantum groups are related to unrelated, at first sight, 
areas of both physics and mathematics 
(Refs. \cite{2dim}--\cite{ky1} and 
references therein). 
For applications of quantum groups, as in the non--deformed
cases, we often need their explicit representations. 
However, despite of remarkable results in this direction 
the problem of
investigating and constructing explicit representations of 
quantum groups, especially those for quantum superalgebras, is 
still far from being satisfactorily solved. Even in the case of 
one--parametric quantum superalgebras, explicit representations 
are mainly known for quantum Lie superalgebras of lower ranks and 
of particular types like $U_{q}[osp(1/2)]$, $U_{q}[gl(1/n)]$ 
(Refs. \cite{ky1,cpt,pt}), while for higher rank quantum Lie 
superalgebras \cite{man3}--\cite{ky3}, besides some q--oscillator
representations which are most popular among those constructed, we 
do not know so much about other representations, in particular, the 
finite--dimensional ones which in many cases are related to 
trigonometric solutions of the quantum Yang--Baxter equations 
\cite{collect,2dim}. Some general aspects and the module 
structure finite--dimensional representations of the quantum 
superalgebras $U_{q}[gl(m/n)]$ are considered in Ref. \cite{zhang} 
but, unfortunately, their explicit construction is still absent. 
So far explicit finite--dimensional irreducible representations 
are all known and classified only for those $U_{q}[gl(m/n)]$ with 
both $m$ and $n\leq 2$ (see Refs. \cite{ky1,ky2,ky3}).\\
     
 What about multi--parametric deformations (first considered
in \cite{man}), this area is even less covered and results are 
much poorer. Some kinds of two--parametric deformations have 
been considered by several authors from different points of 
view (see Refs.\cite{dobrev,ky4} and references therein) but, 
to our knowledge, explicit representations are known and/or 
classified in a few lower rank cases such as 
$U_{p,q}[sl(2/1)]$ and $U_{p,q}[gl(2/1)]$ only 
\cite{ky4,zhang2}. The latter two--parametric quantum 
superalgebra $U_{p,q}[gl(2/1)]$ was consistently defined and 
investigated in \cite{ky4} where all its finite--dimensional 
irreducible representations were explicitly constructed and 
classified at generic deformation parameters. This 
$U_{p,q}[gl(2/1)]$, however, is still a small quantum 
superalgebra which can be defined without the so--called 
extra--Serre defining relations \cite{extra1,extra2,extra3}
representing additional constraints on odd Chevalley generators 
in higher rank cases. Now, in order to include the extra--Serre 
relations on examination we consider a bigger two--parametric 
quantum superalgebra, namely $U_{p,q}[gl(2/2)]$, and its 
representations. This quantum superalgebra $U_{p,q}[gl(2/2)]$ 
resembles to the one--parametric quantum superalgebra 
$U_{\sqrt{pq}}[gl(2/2)]$ but can not be
identified with the latter. Here we suppose $p\neq q$, 
otherwise we 
should return to the case of $U_{q}[gl(2/2)]$ investigated already  
in \cite{ky2,ky3}. Another our motivation is that already in the 
non--deformed case, the superalgebras $gl(n/n)$, especially, their 
subalgebras $sl(n/n)$ and $psl(n/n)$, have special properties (in 
comparison with other $gl(m/n)$, $m\neq n$) and, therefore, attract 
interest \cite{sigma}. Additionally, structures of two--parameter 
deformations considered in \cite{ky4} and
here are, of course, richer than those of one--parameter 
deformations. Every deformation parameter can be independently 
chosen to take a separate generic value (including zero) or 
to be a root of unity.\\
     
  Combining the advantages of the previously developed methods
for $U_{q}[gl(2/2)]$ and $U_{p,q}[gl(1/2)]$ (see Refs. 
\cite{ky2,ky3,ky4}) we can construct all 
finite--dimensional representations of the two--parametric 
quantum Lie superalgebra $U_{p,q}[gl(2/2)]$.
In the frame--work of this paper 
we consider representations at generic 
$p$ and $q$ only (i.e., $p$ and $q$ are not roots of unity), 
while representations at roots of unity are a subject
of later separate investigations. In comparison with previous 
papers \cite{ky2,ky4}, the approach here is somewhat modified 
because of some specific features arising in the present case
but the main steps in the construction procedure remain the same. 
Following this approach we can directly construct explicit 
representations of
the quantum superalgebra $U_{p,q}[gl(2/2)]$ induced from some 
(usually, irreducible) finite--dimensional representations of 
the even subalgebra $U_{p,q}[gl(2)\oplus gl(2)]$ which itself
is a quantum algebra. Since the latter is a stability subalgebra 
of $U_{p,q}[gl(2/2)]$ we expect the representations
of $U_{p,q}[gl(2/2)]$ constructed are decomposed into 
finite--dimensional irreducible representations of 
$U_{p,q}[gl(2)\oplus gl(2)]$. For a clear description of this 
decomposition we shall introduce a $U_{p,q}[gl(2/2)]$--basis 
(i.e., a basis within a $U_{p,q}[gl(2/2)]$--module or briefly 
a basis of $U_{p,q}[gl(2/2)]$) which will be convenient for us 
in investigating the module structure. This basis (see (4.26)) 
can be expressed in terms of some basis of the even subalgebra
$U_{p,q}[gl(2)\oplus gl(2)]$ which in turn represents a (tensor) 
product between two $U_{p,q}[gl(2)]$--bases referred to as the 
left and the right ones. As is shown in \cite{ky4}, the 
Gel'fand--Zetlin (GZ) patterns can serve again as a basis of 
finite--dimensional representations of $U_{p,q}[gl(2)]$. Thus, 
finite--dimensional representations of
$U_{p,q}[gl(2)\oplus gl(2)]$ are realized in tensor products 
of two such $U_{p,q}[gl(2)]$ GZ bases.
For generic $p$ and $q$, the finite--dimensional 
$U_{p,q}[gl(2/2)]$--modules constructed have similar 
structures to those of $U_{q}[gl(2/2)]$ investigated in 
\cite{ky2,ky3} and to those of $gl(2/2)$ investigated 
in \cite{ky6}. However, finite--dimensional representations of 
$U_{p,q}[gl(2/2)]$ at generic deformation parameters are 
not simply trivial deformations from those of $gl(2/2)$ 
that is the former can not be obtained from the latter by 
putting quantum deformation brackets in appropriate places, 
unlike many cases of one--parametric deformations. When
one or both of $p$ and $q$ are roots of unity the structures 
of $U_{p,q}[gl(2/2)]$--modules are drastically different but 
we hope that the present method for construction of 
finite--dimensional representations of $U_{p,q}[gl(2/2)]$
at generic deformation parameters can be extended on its 
finite--dimensional representations at roots of unity.\\
     
  This paper is organized as follows.
The two--parametric quantum superalgebra $U_{p,q}[gl(2/2)]$
is consistently defined in section 2 where we also describe how 
to construct its representations induced from representations of 
the stability subalgebra $U_{p,q}[gl(2)\oplus gl(2)]$. Section 3 
is devoted to constructing finite--dimensional representations of 
$U_{p,q}[gl(2/2)]$. Finally, some comments and conclusions are 
made in section 4.\\
     
 Let us list some abbreviations and notations used
throughout the paper:
\begin{tabbing}
\=12345678\=$V_{l}\otimes V_{r}$ - \=tensor product between two 
linear spaces
$V_{l}$ and $V_{r}$\= or a tensor product\=\kill
\>\>fidirmod(s) : finite-dimensional irreducible module(s),\\[2mm] 
\>\>GZ basis : Gel'fand--Zetlin basis,\\[2mm]
\>\>QGZ basis : Quasi--Gel'fand--Zetlin basis,\\[2mm] 
\>\>lin.env.\{X\} : linear envelope of X,\\[2mm] 
\>\>$p,q$ : the deformation parameters,\\[2mm]
\>\>$[x]\equiv [x]_{p,q}={q^{x}-p^{-x}\over q-p^{-1}}$,~ : 
a $pq$--deformation of a number or an operator $x$,\\[2mm] 
\>\>$V^{p,q}_{l}\otimes V^{p,q}_{r}$ : 
a tensor product between two linear spaces
$V^{p,q}_{l}$ and $V^{p,q}_{r}$\\
\>\>\>~~~~~or a tensor product between a 
$U_{p,q}[gl(2)_{l}]$--module 
$V^{p,q}_{l}$ and\\
\>\>\>~~~~~a $U_{p,q}[gl(2)_{r}]$--module 
$V^{p,q}_{r}$,\\[2mm] 
\>\>$T^{p,q}\odot V^{p,q}_{0}$ : a tensor product between 
two $U_{p,q}[gl(2)\oplus
gl(2)]$--modules\\
\>\>\>~~~~~$T^{p,q}$ and $V^{p,q}_{0}$,\\[2mm]
\>\>$[E,F\}$ : supercommutator between $E$ and $F$,\\[2mm] 
\>\>$[E,F]_{r}\equiv EF-rFE$ : an r--deformed commutator between 
$E$ and $F$,\\[2mm]
\end{tabbing} We hope the notations $[x]\equiv [x]_{p,q}$ for quantum 
deformations, $[m]$ for highest weights (signatures) in
(quasi--) GZ bases $(m)$ and [ , ] for commutators do not confuse 
the reader.
\\[7mm]
{\large {\bf II. The quantum superalgebra $U_{p,q}[gl(2/2)]$}}\\
     
  The quantum superalgebra $U_{p,q}\equiv U_{p,q}[gl(2/2)]$
as a two--parametric deformation of the universal enveloping 
algebra $U[gl(2/2)]$ of the Lie superalgebra $gl(2/2)$ is 
generated by the operators $L_k$, $E_{12}$, $E_{23}$, $E_{34}$, 
$E_{21}$, $E_{32}$, $E_{43}$ and $E_{ii}$ ($1\leq i \leq 4$)
called again Cartan--Chevalley generators and satisfying \cite{ky5} 
\begin{tabbing}
\=123456$[E_{12},E_{34}]$\=$[E_{ii},E_{jj}]$12345678912\= =1\= 0,1234 
\=$[E_{ii},E_{j,j+1}]$\= =
\=$(\delta_{ij}-\delta_{i,j+1})1234$\=\kill
     
{}~~~~a) {\it the super--commutation relations} 
($1\leq i,i+1,j,j+1\leq
4$):\\[2mm]
\>\>$[E_{ii},E_{jj}]$\> = \>0,\>\>\>\>(2.1a) \\[1mm] 
\>\>$[E_{ii},E_{j,j+1}]$\>=\>$(\delta_{ij}-\delta_{i,j+1})E_{j,j+1},$ 
\>\>\>\>(2.1b)\\[1mm]
\>\>$[E_{ii},E_{j+1,j}]$\>=\>$(\delta_{i,j+1}-\delta_{ij})E_{j+1,j}$, 
\>\>\>\>(2.1c)\\[1mm]
     
\>\>[even generator, $L_{k}$]\>=\>0,~ $k=1,2,3$,~
\>\>\>\>(2.1d)\\[1mm]
     
\>\>$[E_{i,i+1},E_{j+1,j}\}$\>=\>$\delta_{ij} 
\left({q\over p}\right)^{L_{i}-h_{i}(1+\delta_{i2})/2} 
[h_{i}]$,
\>\>\>\>(2.1e)\\[4mm]
with
$h_{i}=(E_{ii}-{d_{i+1}\over d_{i}}E_{i+1,i+1})$, 
$L_{1}=L_{l}, L_{2}=0, L_{3}=L_{r}$ and 
$d_{1}=d_{2}=-d_{3}=-d_{4}$\\
$=1$,
\end{tabbing}
     
\begin{tabbing}
\=123456781$[E_{12},E_{34}]$=\=$[E_{ii},E_{jj}]$1234\= =1\= 0,1234 
\=$[E_{ii},E_{j,j+1}]$\= =
\=$(\delta_{ij}-\delta_{i,j+1})E_{j,j+1}$12\=\kill 
{}~~~~b) {\it the Serre--relations}:\\[2mm]
 \>\>~~$[E_{12},E_{34}]$\>=\>$[E_{21},E_{43}]$\>~~~~~~~=~~0,
\>\>\>(2.2a)\\[1mm]
\>\>~~~~~~~$E_{23}^{2}$\>=\>~~~~~$E_{32}^{2}$\>~~~~~~~=~~0, 
\>\>\>(2.2b)\\[1mm]
     
\>~~~~~~~~~~
$[E_{12},E_{13}]_{p}$~~=\> ~~$[E_{21},E_{31}]_{q}$\>=~~ 
\>$[E_{24},E_{34}]_{q}$\>~~~~~~~=~~$[E_{42},E_{43}]_{p}=~~0$,
                           \>\>\>(2.2c)
\end{tabbing}
and\\[2mm]
{}~~~~c) {\it the extra--Serre relations}:\\[2mm] 
$$\{E_{13},E_{24}\}=0,\eqno(2.3a)$$
$$\{E_{31},E_{42}\}=0.\eqno(2.3b)$$
Here, the operators
\\[2mm]
$$E_{13}~:=~[E_{12},E_{23}]_{q^{-1}},\eqno(2.4a)$$ 
$$E_{24}~:=~[E_{23},E_{34}]_{p^{-1}},\eqno(2.4b)$$ 
$$~~~E_{31}~:=~-[E_{21},E_{32}]_{p^{-1}},\eqno(2.4c)$$ 
$$~~~E_{42}~:=~-[E_{32},E_{43}]_{q^{-1}}\eqno(2.4d)$$ 
and the operators composed in the following way 
\begin{tabbing}
\=123456789123456\=$E_{41}$~\=:=1
\=$[E_{21},[E_{32},E_{43}]_{q^{-2}}]_{q^{-2}}$
\=$\equiv [E_{21},E_{42}]_{q^{-1}}$1234567891234~\=\kill 
\>\>$E_{14}$~\>:=\>$[E_{12},[E_{23},E_{34}]_{p^{-1}}]_{q^{-1}}$ 
\>$\equiv ~[E_{12},E_{24}]_{q^{-1}}$,\>(2.5a)\\[1mm] 
\>\>$E_{41}$\>:=\>$[E_{21},[E_{32},E_{43}]_{q^{-1}}]_{p^{-1}}$ 
\>$\equiv ~-[E_{21},E_{42}]_{p^{-1}}$\>(2.5b)
\end{tabbing}
are defined as new generators which, like $E_{23}$ and 
$E_{32}$, are all odd and have vanishing squares. 
These generators $E_{ij}$, $1\leq i,j\leq 4$, are
two--parametric deformation analogues ($pq$--analogues) of the 
Weyl generators $e_{ij}$, $1\leq i,j\leq 4$, of the superalgebra 
$gl(2/2)$ whose universal enveloping algebra $U[gl(2/2)]$ is a 
classical limit of $U_{p,q}[gl(2/2)]$ when $p,q\rightarrow 1$. 
The so--called maximal--spin operators $L_k$ are constants 
within a $U_{p,q}[gl(2)]$--fidirmod and are different for 
different $U_{p,q}[gl(2)]$--fidirmods. Therefore, commutators 
between these operators with the odd generators intertwining 
$U_{p,q}[gl(2)]$--fidirmods take concrete forms on concrete 
basis vectors. Other commutation relations between $E_{ij}$ 
follow from the relations (2.1)--(2.3) and the definitions 
(2.4) and (2.5).\\
     
 The subalgebra $U_{p,q}[gl(2/2)_{0}]~ (\subset
U_{p,q}[gl(2/2)]_{0}\subset U_{p,q}[gl(2/2)])$ is even and 
isomorphic to $U_{p,q}[gl(2)\oplus gl(2)]\equiv U_{p,q}[gl(2)] 
\oplus U_{p,q}[gl(2)]$ which is completely defined by $L_1$, 
$L_3$, $E_{12}$, $E_{34}$, $E_{21}$, $E_{43}$ and $E_{ii}$, 
$1\leq i\leq 4$,
$$U_{q}[gl(2/2)_{0}]~=~ {\normalsize lin.env.} \{L_1, L_3, E_{ij}\|~ 
i,j=1,2~~
{\normalsize and}~~i,j=3,4\}.\eqno(2.6)$$ 
In order to distinguish two components $U_{p,q}[gl(2)]$ 
of $U_{p,q}[gl(2/2)_{0}]$ we set 
$$left~~U_{p,q}[gl(2)]\equiv U_{p,q}[gl(2)_{l}]:=
{\normalsize lin.env.}\{L_1, E_{ij}\|~i,j=1,2\},\eqno(2.7)$$ 
$$right~U_{p,q}[gl(2)]\equiv U_{p,q}[gl(2)_{r}]:= 
{\normalsize lin.env.}\{L_3, E_{ij}\|~ i,j=3,4\}, 
\eqno(2.8)$$
that is
$$U_{p,q}[gl(2/2)_{0}]~=~U_{p,q}[gl(2)_{l}\oplus gl(2)_{r}]. 
\eqno(2.9)$$
     
 Looking at the relations (2.1)--(2.3) we see that every of
the odd spaces $A_{+}$ and $A_{-}$ spanned on the positive 
and negative odd roots (generators) $E_{ij}$ and $E_{ji}$, 
$1\leq i\leq 2 <j\leq 4$, respectively
$$A_{+}= {\normalsize lin.env.}\{E_{14},E_{13}, 
E_{24}, E_{23}\},\eqno(2.10)$$
$$A_{-}= {\normalsize lin.env.}\{E_{41},E_{31}, E_{42},E_{32} \}, 
\eqno(2.11)$$
is a representation space of the even subalgebra 
$U_{p,q}[gl(2/2)_{0}]$
which, as seen from (2.1)--(2.2), is a stability subalgebra 
of $U_{p,q}[gl(2/2)]$.
Therefore, we can construct representations of 
$U_{p,q}[gl(2/2)]$ induced from some (finite--dimensional 
irreducible, for example) representations of 
$U_{p,q}[gl(2/2)_{0}]$ which are realized
in some representation spaces (modules) $V^{p,q}_{0}$ 
representing tensor products of
$U_{p,q}[gl(2)_l]$--modules $V^{p,q}_{0,l}$ and 
$U_{p,q}[gl(2)_r]$--modules $V^{p,q}_{0,r}$ 
$$V_{0}^{p,q}(\Lambda)=V_{0,l}^{p,q}(\Lambda_{l})\otimes 
V_{0,r}^{p,q}(\Lambda_{r}),\eqno(2.12)$$
where  $\Lambda$'s are some signatures (such as highest 
weights, respectively) characterizing the modules (highest 
weight modules, respectively). Here $\Lambda_l$ and 
$\Lambda_r$ are referred to as the left and the right 
components of $\Lambda$, respectively
     
$$\Lambda=[\Lambda_{l},\Lambda_{r}].\eqno(2.13)$$ 
     
  If we demand
$$E_{23}V_{0}^{p,q}(\Lambda)=0 \eqno(2.14)$$ hence 
$$U_{p,q}(A_{+})V_{0}^{p,q}=0,\eqno(2.15)$$ we 
turn the $U_{p,q}[gl(2/2)_{0}]$--module 
$V^{p,q}_{0}$ into a $U_{p,q}(B)$--module where 
$$B=A_{+}\oplus gl(2)\oplus gl(2).\eqno(2.16)$$
The $U_{p,q}[gl(2/2)]$--module $W^{p,q}$ induced from 
the $U_{p,q}[gl(2/2)_{0}]$--module $V^{p,q}_{0}$ is the 
factor--space
$$W^{p,q}=W^{p,q}(\Lambda)=[U_{p,q}\otimes V_{0}^{p,q} 
(\Lambda)]/I^{p,q}(\Lambda)\eqno(2.17)$$
which, of course, depends on $\Lambda$, where 
$$U_{p,q}\equiv U_{p,q}[gl(2/2)], \eqno(2.18)$$ 
while $I^{p,q}$ is the subspace
$$I^{p,q}={\normalsize lin.env.}
\{ub\otimes v-u\otimes bv\| u\in U_{p,q},
b\in U_{p,q}(B)\subset U_{p,q}, v\in V_{0}^{p,q}\}. 
\eqno(2.19)$$
   Using the commutation relations (2.1)--(2.3) and the
definitions (2.4) and (2.5) we can prove the following 
analogue of the Poincar\'e--Birkhoff--Witt theorem\\[4mm] 
{\bf Proposition 1}: {\it The quantum deformation
$U_{p,q} := U_{p,q}[gl(2/2)]$ is spanned on all possible 
linear combinations of the elements
$$g=(E_{23})^{\eta_{1}}(E_{24})^{\eta_{2}} 
(E_{13})^{\eta_{3}}(E_{14})^{\eta_{4}}(E_{41})^{\theta_{1}} 
(E_{31})^{\theta_{2}}(E_{42})^{\theta_{3}}(E_{32})^ 
{\theta_{4}}g_{0}, \eqno(2.20)$$
or equivalently
$$g=(E_{41})^{\theta_{1}}(E_{31})^{\theta_{2}} 
(E_{42})^{\theta_{3}}(E_{32})^{\theta_{4}}b, \eqno(2.21)$$
where $g_{0}\in U_{p,q}[gl(2/2)_{0}]$, $b\in U_{p,q}(B)$ 
and $\eta_{i}$, $\theta_{i} = 0,1$}.\\
     
   Any vector $w$ from the module $W^{p,q}$ can be represented 
as
$$w=u\otimes v,~~~~ u\in U_{p,q},~~~~ 
v\in V_{0}^{p,q}.\eqno(2.22)$$ 
Then $W^{p,q}$ is a $U_{p,q}[gl(2/2)]$--module in the sense 
$$gw\equiv g(u\otimes v)=gu\otimes v\in W^{p,q},\eqno(2.23)$$
for $g,~u\in U_{p,q}$, $w\in W^{p,q}$ and $v\in V_{0}^{p,q}$. 
Taking into account the fact that $V_{0}^{p,q}(\Lambda)$ is a 
$U_{p,q}(B)$--module we have
$$W^{q}(\Lambda)={\normalsize lin.env.}
\{(E_{41})^{\theta_{1}}(E_{31})^{\theta_{2}} 
(E_{42})^{\theta_{3}}(E_{32})^{\theta_{4}}\otimes v\|~v\in 
V_{0}^{p,q},~\theta_{1},...,\theta_{4}=0,1\}.\eqno(2.24)$$ 
Consequently, a basis of $W^{p,q}$ can be constituted by
taking all the vectors of the form
$$\left |\theta_{1}, \theta_{2}, \theta_{3}, \theta_{4}; 
(\lambda)\right > :=
(E_{41})^{\theta_{1}}
(E_{31})^{\theta_{2}}(E_{42})^{\theta_{3}}(E_{32})^ 
{\theta_{4}}\otimes (\lambda),
~~ \theta_{i}=0,1, \eqno(2.25)$$ where $(\lambda)$ 
is a (GZ, for example) basis of $V_{0}^{p,q}\equiv 
V_{0}^{p,q}(\Lambda)$. We refer to this basis of 
$W^{p,q}$ as the induced $U_{p,q}[gl(2/2)]$--basis
(or simply, the induced basis) in order to distinguish 
it from another $U_{p,q}[gl(2/2)]$--basis introduced
later and called a reduced basis which is more convenient 
for investigating the module structure of $W^{p,q}$.
It is obvious that if the module $V_{0}^{p,q}$ is 
finite--dimensional so is the module $W^{p,q}$. 
Finite--dimensional representations of $U_{p,q}[gl(2/2)]$ 
are namely the subject of the next section.\\[7mm]
{\large {\bf III. Finite--dimensional representations of 
$U_{p,q}[gl(2/2)]$}}\\
     
   In this section we consider
finite--dimensional representations of $U_{p,q}[gl(2/2)]$ 
induced from irreducible finite--dimensional representations 
of $U_{p,q}[gl(2/2)_{0}]$. We firstly construct the bases of 
the module $W^{q}$ and then find the explicit matrix elements
for the finite--dimensional representations of $U_{p,q}[gl(2/2)]$. 
\\
     
  We can shown that the GZ patterns
$$
\left[
\begin{array}{c}
                        m_{12}~~~m_{22}\\ m_{11}
\end{array}
\right]
\equiv
\left[
\begin{array}{c}
                             [m]\\ m_{11}
\end{array}
\right]
\eqno(3.1)$$
where $m_{ij}$ are complex numbers such that $m_{12}-m_{11}\in 
{\bf Z}_{+}$ and $m_{11}-m_{22}\in {\bf Z}_{+}$, can serve as 
a basis of a $U_{p,q}[gl(2)]$--fidirmod. Indeed, 
finite--dimesional representations of $U_{p,q}[gl(2)]$
are high weight and if the operators $L$ and $E_{ij}$, 
$i,j=1,2,$ are difined on the basis (3.1) as follows 
\begin{eqnarray*}
~~~~~L\left[
\begin{array}{c}
                        m_{12}~~~m_{22}\\ m_{11}
\end{array}
\right]
& = &{1\over 2}(l_{12}-l_{22}-1)\left[
\begin{array}{c}
                        m_{12}~~~m_{22}\\ m_{11}
\end{array}
\right],\\[2mm]
E_{11}\left[
\begin{array}{c}
                        m_{12}~~~m_{22}\\ m_{11}
\end{array}
\right]
& = &(l_{11}+1)\left[
\begin{array}{c}
                        m_{12}~~~m_{22}\\ m_{11}
\end{array}
\right],\\[2mm]
E_{22}\left[
\begin{array}{c}
                        m_{12}~~~m_{22}\\ m_{11}
\end{array}
\right]
& = &(l_{12}+l_{22}-l_{11}+2)\left[
\begin{array}{c}
                        m_{12}~~~m_{22}\\ m_{11}
\end{array}
\right],\\[2mm]
E_{12}\left[
\begin{array}{c}
                        m_{12}~~~m_{22}\\ m_{11}
\end{array}
\right]
& = &\left([l_{12}-l_{11}][l_{11}-l_{22}]\right)^{1/2} 
\left[
\begin{array}{c}
                        m_{12}~~~m_{22}\\ m_{11}+1
\end{array}
\right],\\[2mm]
E_{21}\left[
\begin{array}{c}
                        m_{12}~~~m_{22}\\ m_{11}
\end{array}
\right]
& = &\left([l_{12}-l_{11}+1][l_{11}-l_{22}-1]\right)^{1/2} 
\left[
\begin{array}{c}
                        m_{12}~~~m_{22}\\ m_{11}-1
\end{array}
\right],~~~~~~~(3.2)
\end{eqnarray*}
they really satisfy commutation relations of $U_{p,q}[gl(2)]$ 
given in (2.1). Here the notation
     
$$l_{ij}=m_{ij}-i, ~~ i=1,2\eqno(3.3a)$$ 
and later also the notation
     
$$l'_{ij}=m_{ij}-i+2, ~~i=3,4,\eqno(3.3b)$$
are used. Since the $U_{p,q}[gl(2/2)_{0}]$--fidirmod 
$V_{0}^{p,q}$ is decomposed into a
$U_{p,q}[gl(2)_{l}]$--fidirmod $V_{0,l}^{p,q}$ and a 
$U_{p,q}[gl(2)_{r}]$--fidirmod $V_{0,r}^{p,q}$ via the 
tensor product
$$V_{0}^{p,q}=V_{0,l}^{p,q}\otimes V_{0,r}^{p,q},\eqno(3.4)$$ 
its basis, therefore, is a tensor product
$$
\left[
\begin{array}{c}
                        m_{13}~~~m_{23}\\ m_{11}
\end{array}
\right]
\otimes
\left[
\begin{array}{c}
                        m_{33}~~~m_{43}\\ m_{31}
\end{array}
\right]
\equiv
\left[
\begin{array}{c}
                            [m]_{l}\\ m_{11}
\end{array}
\right]
\otimes
\left[
\begin{array}{c}
                            [m]_{r}\\ m_{31}
\end{array}
\right]
\equiv
(m)_{l}\otimes (m)_{r}
\equiv
(m)
\eqno(3.5)$$
between a GZ basis of $V_{0,l}^{p,q}$ spanned on the vectors 
$(m)_{l}$ and a GZ basis  of $V_{0,r}^{p,q}$ spanned on the 
vectors $(m)_{r}$.  Following the approach of Ref. \cite{ky2} 
(see aslo Ref.\cite{ky3,ky6}) and keeping the notations used 
there, we can represent the basis (3.5) of $V_{0}^{p,q}$ in 
the form
$$\left[
\begin{array}{lcr}
     
\begin{array}{c}
                        m_{13}~~~m_{23}\\ m_{11}
\end{array}
;
\begin{array}{c}
                        m_{33}~~~m_{43}\\ m_{31}
\end{array}
\end{array}
\right]
\equiv
\left[
\begin{array}{lcr}
     
\begin{array}{c}
                            [m]_{l}\\ m_{11}
\end{array}
;
\begin{array}{c}
                            [m]_{r}\\ m_{31}
\end{array}
\end{array}
\right]
\equiv
(m)
\eqno(3.6)$$
 Then, the signature $\Lambda$, which now is the highest weight,
is given by the first row
$[m_{13},m_{23},m_{33},m_{43}]\equiv [[m]_{l},[m]_{r}]\equiv [m]$ 
common for all the basis vectors (3.6) of $V_{0}^{p,q}$: 
$$V_{0}^{p,q}\equiv
V_{0}^{p,q}(\Lambda)=V_{0}^{p,q}([m])= V_{0,l}^{p,q}([m]_{l})\otimes 
V_{0,r}^{p,q}([m]_{r}).\eqno(3.7)$$ The explicit action of 
$U_{p,q}[gl(2/2)_{0}]$ on $V_{0}^{p,q}([m])$ follows directly from 
(3.2) and : $$g_{0}(m)=g_{0,l}(m)_{l}\otimes (m)_{r} + 
(m)_{l}\otimes g_{0,r}(m)_{r}\eqno(3.8)$$ for $g_{0}\equiv 
g_{0,l}\oplus g_{0,r}\in U_{q}[gl(2/2)_{0}]$ and $(m)\in 
V_{0}^{q}([m])$.\\
     
   The basis vector with $ m_{11}= m_{13}$ and 
$ m_{31}= m_{33}$
     
$$
\left[
\begin{array}{c}
                        m_{13}~~~m_{23}\\ m_{13}
\end{array}
;
\begin{array}{c}
                        m_{33}~~~m_{43}\\ m_{33}
\end{array}
\right]
\equiv
\left[
\begin{array}{c}
                            [m]_{l}\\ m_{13}
\end{array}
;
\begin{array}{c}
                            [m]_{r}\\ m_{33}
\end{array}
\right]
\equiv (M)\eqno(3.9)$$
satisfying the conditions
\begin{tabbing}
\=123456789123456789123456\=$E_{ii}(M)$~\==~\=$m_{i3}(M)$ 
~~~~$i=1,2,3,4$123456789123.\=\kill
\>\>$E_{ii}(M)$\>=\>$m_{i3}(M)$,~~~~ 
$i=1,2,3,4$,\\[1mm]
\>\>$E_{12}(M)$\>=\>$E_{34}(M)~=~0$\>~~~(3.10) 
\end{tabbing}
is the highest weight vector in $V_{0}^{p,q}([m])$. 
Therefore,
as in the classical case ($p=q=1$) \cite{ky6} and in
the case of one--parametric deformation ($p=q$) \cite{ky2} 
the highest weight $[m]$ is nothing but an ordered set of
the eigen values of the Cartan generators $E_{ii}$ on the highest 
weight vector $(M)$. The latter is also highest weight vector in 
$W^{p,q}([m])$ because of the condition (2.14). All other, i.e. 
lower weight, basis vectors of $V_{0}^{p,q}$ can be obtained from 
the highest weight vector $(M)$ through acting on the latter by 
monomials of the lowering generators $E_{21}$
and $E_{43}$ in definite powers:
\begin{eqnarray*}
~~~~~~~~(m)&=&\left({[m_{11}-m_{23}]![m_{31}-m_{43}]!\over 
[m_{13}-m_{23}]![m_{13}-m_{11}]!
[m_{33}-m_{43}]![m_{33}-m_{31}]!}\right)^{1/2}\\[5mm]&&\times 
(E_{21})^{m_{13}-m_{11}}(E_{43})^{m_{33}-m_{31}}(M), 
~~~~~~~~~~~~~~~~~~~~~~~~~~~~~~~~~~~~~~(3.11)
\end{eqnarray*}
where $[n]$'s stand for $${q^{n}-p^{-n}\over 
q-p^{-1}}\equiv [n]_{p,q}\equiv [n],\eqno(3.12)$$ 
while
$$[n]!=[1][2]...[n-1][n].\eqno(3.13)$$
Therefore, the induced basis (2.25) of $W^{p,q}(\Lambda)=W^{p,q}([m])$ 
now takes the form
$$\left|\theta_{1}\theta_{2},\theta_{3},\theta_{4};(m)\right>~:= 
(E_{41})^{\theta_{1}}(E_{31})^{\theta_{2}} 
(E_{42})^{\theta_{3}}(E_{32})^{\theta_{4}}\otimes 
(m).\eqno(3.14)$$
The subspace $T^{p,q}$ consisting of 
$$\left|\theta_{1},\theta_{2},\theta_{3},\theta_{4}\right>~:= 
(E_{41})^{\theta_{1}}(E_{31})^{\theta_{2}} 
(E_{42})^{\theta_{3}}(E_{32})^{\theta_{4}}\eqno(3.15)$$ can be 
considered as a $U_{p,q}[gl(2/2)_{0}]$--adjoint module which is 
16--dimensional when all the $\theta_i$ ($i=1,2,3,4$) take all 
two possible values 0 and 1, that is $\sum_{i=1}^{4}\theta_i$ 
runs all over the range from 0 to 4.
Thus $W^{p,q}([m])$ being a tensor product between 
two $U_{p,q}[gl(2/2)_{0}]$--modules:
$$W^{p,q}([m])=T^{p,q}\odot V_{0}^{p,q}([m]),\eqno(2.24')$$ 
is, in general, a reducible $U_{p,q}[gl(2/2)_{0}]$--module 
and are decomposed into irreducible
$U_{p,q}[gl(2/2)_{0}]$--submodules.
We arrive at the next assertion\\[4mm]
{\bf Proposition 2}: {\it The induced $U_{p,q}[gl(2/2)]$--module 
$W^{p,q}$ is the linear span
$$W^{p,q}([m])={\normalsize
lin.env.}\{(E_{41})^{\theta_{1}}(E_{31})^{\theta_{2}} 
(E_{42})^{\theta_{3}}(E_{32})^{\theta_{4}}\otimes v\|v\in 
V_{0}^{p,q}([m]),~\theta_{i}=0,1\},\eqno(2.24'')$$
which is decomposed into a direct sum of (sixteen, at most) 
$U_{p,q}[gl(2/2)_{0}]$--fidirmods $V_{k}^{p,q}([m]_{k})$: 
$$W^{p,q}([m])=\bigoplus_{k=0}^{15}V_{k}^{p,q}([m]_{k}). 
\eqno(3.16)$$
where $[m]_{k}$ are signatures of
$V_{k}^{p,q}\equiv V_{k}^{p,q}([m]_{k})$}.\\ 

  Here, we call $[m]_{k}\equiv [m_{12},m_{22},m_{32},m_{42}]_{k}$
the local highest weights of the submodules $V_{k}^{p,q}$ in their 
GZ bases denoted now as
$$
\left[
\begin{array}{lcr}
     
\begin{array}{c}
                        m_{12}~~~m_{22}\\ m_{11}
\end{array}
;
\begin{array}{c}
                        m_{32}~~~m_{42}\\ m_{31}
\end{array}
\end{array}
\right]_{k}
\equiv
(m)_{k}.
\eqno(3.17)$$
The highest weight $[m]_{0}\equiv [m]$ of $V_{0}^{p,q}$ being also 
the highest weight of $W^{p,q}$ is referred to as the global highest 
weight. We call $[m]_{k}$, $k\neq 0$, the local highest weights
in the sense that they characterize the submodules 
$V_{k}^{p,q}\subset W^{p,q}$ as
$U_{p,q}[gl(2/2)_{0}]$--fidirmods only,
while the global highest weight $[m]$ characterizes the 
$U_{p,q}[gl(2/2)]$--module $W^{p,q}$ as the whole. In the same 
way we define the local highest weight vectors $(M)_{k}$ in 
$V_{k}^{p,q}$ as those $(m)_{k}$ satisfying the conditions (cf. 
(3.10))
\begin{eqnarray*}
~~~~~~~~~~~~~~~~~~~~~~~~~~~~E_{ii}(M)_{k}&=&m_{i2}(M)_{k}, 
~~~~i=1,2,3,4,\\[1mm]
                     E_{12}(M)_{k}&=&E_{34}(M)_{k}~=~0.
~~~~~~~~~~~~~~~~~~~~~~~~~~~~~~~~~~~(3.18) 
\end{eqnarray*}
The highest weight vector $(M)$ of $V_{0}^{p,q}$ is also the global 
highest weight vector in $W^{p,q}$  for which the condition (see 
(2.14))
$$E_{23}(M)=0\eqno(3.19)$$
and  the conditions (3.18)
simultaneously hold.\\
     
  Let us denote by $\Gamma_{k}^{p,q}$ the basis system  spanned
on the basis vectors $(m)_{k}$ (3.17) in each $V_{k}^{p,q}([m])$. 
For a basis of  $W^{p,q}$ we can choose the union 
$\Gamma^{p,q}=\bigcup_{k=0}^{15}\Gamma_{k}^{p,q}$ of all the bases 
$\Gamma_{k}^{p,q}$, namely, a basis vector of $W^{p,q}$ has to be 
identified with one of the vectors $(m)_{k}$,
$0 \leq k \leq15$. The basis $\Gamma^{p,q}$ is referred to as the 
$U_{p,q}[gl(2/2)]$--reduced basis or simply, the reduced basis. 
It is clear that every basis 
$\Gamma_{k}^{p,q}=\Gamma_{k}([m]_{k})^{p,q}$ is labelled by a 
local highest weight $[m]_{k}$, while the basis 
$\Gamma^{p,q}=\Gamma^{p,q}([m])$ is labelled by the global 
highest weight $[m]$. Going ahead, we modify the notation (3.17) 
for the basis vectors in $\Gamma^{p,q}$ as follows (cf. (4.26) 
in Ref. \cite{ky2})
$$
\left[
\begin{array}{lccc}
     
                        m_{13}& m_{23}&  m_{33}& m_{43} \\
     
                        m_{12}& m_{22}&  m_{32}&  m_{42}\\
m_{11}&   0   &  m_{31}&    0
\end{array}
\right]_{k}
\equiv
\left[
\begin{array}{lcr}
     
\begin{array}{c}
                        m_{12}~~~m_{22}\\ m_{11}
\end{array}
;
\begin{array}{c}
                        m_{32}~~~m_{42}\\ m_{31}
\end{array}
\end{array}
\right]_{k}
\equiv
(m)_{k},
\eqno(3.20)$$
with $k$ running from 0 to 15 as for $k=0$ we must take into 
account $m_{i2}=m_{i3}$, $i=1,2,3,4$, i.e., 
$$ (m)_{0}\equiv (m)=
\left[
\begin{array}{lccc}
     
                        m_{13}& m_{23}&  m_{33}& m_{43} \\
     
                        m_{13}& m_{23}&  m_{33}&  m_{43}\\
m_{11}&   0   &  m_{31}&    0
\end{array}
\right].
\eqno(3.21)
$$ In (3.20) the first row $[m]=[m_{13},m_{23},m_{33},m_{43}]$ 
being the (global) highest weight of $W^{p,q}$ is fixed for all the
vectors in the whole $W^{p,q}$ and characterizes this module itself, 
while the second row is a (local) highest weight of some
submodule $V_{k}^{p,q}$ and tells us that the considered basis 
vector $(m)_{k}$ of $W^{p,q}$ belongs to this submodule in the 
decomposition (3.16) corresponding to the branching rule 
$U_{p,q}[gl(2/2)]\supset U_{p,q}[gl(2/2)_{0}]\supset 
U_{p,q}[gl(1)\otimes gl(1)]$. We refer to
(3.20) as the quasi--Gel'fand--Zetlin (QGZ) basis.\\

   It is easy to  see that the highest vectors $(M)_{k}$  in the
notation (3.20) are
$$ (M)_{k}=
\left[
\begin{array}{lccc}
     
                        m_{13}& m_{23}&  m_{33}& m_{43} \\
     
                        m_{12}& m_{22}&  m_{32}&  m_{42}\\
m_{12}&   0   &  m_{32}&    0
\end{array}
\right]_{k},~~~k=0,1,...15.
\eqno(3.22)$$
The (global) highest weight vector $(M)$ (3.9) is given now by 
$$ (M)=
\left[
\begin{array}{lccc}
     
                        m_{13}& m_{23}&  m_{33}& m_{43} \\
     
                        m_{13}& m_{23}&  m_{33}&  m_{43}\\
m_{13}&   0   & m_{33}& 0 \end{array} \right].  \eqno(3.23)$$ 
A highest weight vector $(M)_{k}$ expressed in 
terms of the induced basis (3.14) has the form of 
a homogeneous polynomial of a definite
degree $\eta$ in negative odd generators
($E_{ij}$, $1\leq j \leq 2 < i \leq 4$) acting on $(m)\in 
V_0^{p,q}([m])$:
$$(M)_k\equiv (M)_{\eta,h}=
\sum_{\theta_i=0,1}
C_{\eta,h}(\theta_{1},\theta_{2},\theta_{3},\theta_{4}) 
\left|\theta_{1},\theta_{2},\theta_{3},\theta_{4};(m)\right> 
\eqno(3.24)$$
with $\eta=\sum_{i=1}^4 \theta_i$ fixed for every 
$(M)_{\eta,h}$, and the coefficients $C_{\eta,h}$ 
determined by solving Eqs. (3.18). Applying (3.11) to 
any $(M)_{\eta,h}$ we find all the basis vectors 
$(m)_{\eta,h}$ of the corresponding fidirmod 
$V_{\eta,h}^{p,q}$ which is a linear space spanned
on homogeneous polynomials of the negative odd generators 
of the same degree $\eta$ since (3.11) does not change
$\eta$. Here we call $\eta$ the level of $V_{\eta,h}^{p,q}$. 
It is easy to see that on the level
$\eta=0$ there is only one fidirmod, namely 
$V_0^{p,q}\equiv V_{0,1}^{p,q}$, while on the next level 
$\eta=1$ there are four fidirmods, say, $V_{1,h}^{p,q},
~h=1,2,3,4$. On the level $\eta=2$ we can find six fidirmods 
$V_{2,h}^{p,q}, ~1\leq h\leq 6$, which are divided into two 
groups ($h=1,2,3$ and $h=4,5,6$) expressed in terms of two 
indenpendent groups of second order monomials of odd 
generators (3.15) acting on $(m)$. For $\eta=3$ the number 
of fidirmods is four, $V_{3,h}^{p,q}, ~h=1,2,3,4$,
and finally on $\eta=4$ we find again only one fidirmod 
$V_{4,1}^{p,q}$. However, this form (3.24) which was used in 
the one--parametric case \cite{ky2} is now inconvenient for 
us here to apply formula (3.11) in order to find all other
(i.e., lower weight) reduced basis vectors. It is so because 
of the presence of the maximal--spin operators $L_i$ which are 
not diagonalized in the induced basis but in the reduced basis
(since an eigenvalue of any $L_i$ is a fixed constant only 
within a $U_{p,q}[gl(2)]$--fidirmod (or fidirmod, for short) 
and changes from fidirmod to fidirmod). Applying (3.11) we 
have to push generators $E_{21}$ and $E_{43}$ to the right 
side until reaching $V_0^{p,q}$ by using commutation 
ralations (2.1)--(2.2) which give rise to $L_i$'s acting 
on the induced basis vectors. 
But it is extremtly difficult to get explicit actions of
$L_i$ on the latter vectors before knowing how they are 
projected on the reduced basis which, however, we are now 
looking for. Instead, we will write down $(M)_k\equiv 
(M)_{\eta,h}$ in a form convenient for applying (3.11) 
which leaves the $\eta$'s unchanged:
\begin{eqnarray*}
~~~~~~(M)_{0}&\equiv &(M)_{0,1}=
a_{0}\left|0,0,0,0;(M)\right> 
\equiv (M),~~~~
a_{0}\equiv
1,\\[2mm] (M)_{1}&\equiv &(M)_{1,1}=
a_{1}\left|0,0,0,1;(M)\right> 
\equiv a_1E_{32}(M),\\[2mm]
(M)_{2}&\equiv &(M)_{1,2}=
a_{2}\left\{{1\over a_1}E_{21}(M)_1 - 
\frac{[2l+1]} {[2l]}
E_{32}E_{21}(M) \right\},\\[2mm]
(M)_{3}&\equiv &(M)_{1,3}=
a_{3}\left\{{1\over a_1}E_{43}(M)_1 - 
\frac{[2l'+1]} {[2l']}
E_{32}E_{43}(M) \right\},\\[2mm]
(M)_{4}&\equiv &(M)_{1,4}=
a_{4}\left\{{1\over a_1}E_{21}E_{43}(M)_1 - 
{1\over a_2} E_{43}(M)_2 -
{1\over a_3} E_{21}(M)_3\right.\\[2mm] 
&&\left.
- \frac{[2l+1][2l'+1]}
{[2l][2l']} E_{32}E_{21}E_{43}(M) \right\},\\[2mm] 
(M)_{5}&\equiv  &(M)_{2,1}=a_{5}\left|0,0,1,1;(M)\right> 
\equiv a_5E_{42}E_{32}(M),\\[4mm]
(M)_{6}&\equiv &(M)_{2,2}=a_{6}\left\{{1\over 
a_5} E_{21}(M)_5 - \frac{[2l+2]} {[2l]}
E_{42}E_{32}E_{21}(M)\right\},\\[2mm]
(M)_{7}&\equiv &(M)_{2,3}=a_{7}\left\{{1\over 
a_5}E_{21}^2(M)_5 - {1 \over a_6}
\frac{[2][2l+1]} {[2l]} E_{21}(M)_6\right.\\[2mm] 
& &\left.- \frac{[2l+1][2l+2]}
{[2l][2l-1]} E_{42}E_{32}E_{21}^2(M) \right\},\\[2mm] 
(M)_{8}&\equiv &(M)_{2,4}=a_{8}\left|0,1,0,1;(M)\right> 
\equiv a_8E_{31}E_{32}(M),\\[2mm]
(M)_{9}&\equiv &(M)_{2,5}=a_{9}\left\{{1\over 
a_8}E_{43}(M)_8 - \frac{[2l'+2]} {[2l']}
E_{31}E_{32}E_{43}(M)\right\},\\[2mm]
(M)_{10}&\equiv & (M)_{2,6}=a_{10}\left\{{1\over 
a_8}E_{43}^2(M)_8 - {1 \over a_9}
\frac{[2][2l'+1]} {[2l']} E_{43}(M)_9\right.\\[2mm] 
& &\left.- \frac{[2l'+1][2l'+2]}
{[2l'][2l'-1]} E_{31}E_{32}E_{43}^2(M) \right\},\\[2mm] 
(M)_{11}&\equiv &(M)_{3,1}=a_{11}\left|0,1,1,1;(M)\right> 
\equiv a_{11}E_{31}E_{42}E_{32}(M),\\[2mm]
(M)_{12}&\equiv &(M)_{3,2}=a_{12}\left\{{1\over 
a_{11}}E_{21}(M)_{11}
- \frac{[2l+1]} {[2l]} E_{31}E_{42}E_{32}E_{21}(M) 
\right\},\\[2mm]
(M)_{13}&\equiv &(M)_{3,3}=a_{13}\left\{{1\over 
a_{11}}E_{43}(M)_{11}
- \frac{[2l'+1]} {[2l']}
E_{31}E_{42}E_{32}E_{43}(M)\right\},\\[2mm] 
(M)_{14}&\equiv &(M)_{3,4}=a_{14}\left\{{1\over 
a_{11}}E_{21}E_{43}(M)_{11} -
{1\over a_{12}} E_{43}(M)_{12} -
{1\over a_{13}} E_{21}(M)_{13}\right.\\[2mm] 
&&\left.
- \frac{[2l+1][2l'+1]}
{[2l][2l']} E_{31}E_{42}E_{32}E_{21}E_{43}(M) 
\right\},\\[2mm]
(M)_{15}&\equiv
&(M)_{4,1}=a_{15}\left|1,1,1,1;(M)\right> 
\equiv a_{15}E_{41}E_{31}E_{42}E_{32}(M), 
~~~~~~~~~~~~(3.25)
\end{eqnarray*}\\ where
$l={1\over 2}(m_{13}-m_{23})$ and $l'=
{1\over 2}(m_{33}-m_{43})$, 
while $a_{k}=a_{k}(p,q)$ are coefficients depending, 
in general, on $p$ and $q$.
Indeed, $(M)_{k}$ given in (3.25) form a set of all 
linear independent vectors satisfying the conditions 
(3.18).
  Looking at (3.25) we  easily
identify the highest weights $[m]_{k}$
\begin{tabbing} \=12345679123456789\= $[m]_{kk}$ \= 
=x \= $[m_{13}-1,m_{23}-1,m_{33}+1,m_{43}+1]$,\=\kill
\>\>
$[m]_{0}$ \> = 
\> $[m_{13}, m_{23}, m_{33}, m_{43}]$,\\[2mm] 
\>\>$[m]_{1}$ \> = \> $[m_{13}, m_{23}-1, m_{33}+1, 
m_{43}]$,\\[2mm] \>\>$[m]_{2}$ \> = \> $[m_{13}-1, 
m_{23}, 
m_{33}+1, m_{43}]$,\\[2mm] \>\>$[m]_{3}$ \> = \> 
$[m_{13}, 
m_{23}-1, m_{33}, m_{43}+1]$,\\[2mm] \>\>
$[m]_{4}$ \> = \> 
$[m_{13}-1, m_{23}, m_{33}, m_{43}+1]$,\\[2mm]
\>\>$[m]_{5}$ \> =
\> $[m_{13},m_{23}-2,m_{33}+1,m_{43}+1]$,\\[2mm] 
\>\>$[m]_{6}$ \>
= \> $[m_{13}-1,m_{23}-1,m_{33}+1,m_{43}+1]_{6}$,\\[2mm] 
\>\>$[m]_{7}$ \> = \>
$[m_{13}-2,m_{23},m_{33}+1,m_{43}+1]$,\\[2mm]
\>\>$[m]_{8}$ \> =\> $[m_{13}-1, m_{23}-1, 
m_{33}+2, m_{43}]$, 
\\[2mm]
\>\>$[m]_{9}$\> = \> $[m_{13}-1, m_{23}-1, m_{33}+1, 
m_{43}+1]_{9}$, 
\\[2mm]
\>\>$[m]_{10}$ \> = \>$[m_{13}-1, m_{23}-1, m_{33}, 
m_{43}+2]$,\\[2mm] 
\>\>$[m]_{11}$ \> =
\> $[m_{13}-1,m_{23}-2,m_{33}+2, 
m_{43}+1]$,\\[2mm] \>\>$[m]_{12}$ 
\> = \> $[m_{13}-2,m_{23}-1,m_{33}+2, 
m_{43}+1]$,\\[2mm] 
\>\>$[m]_{13}$ \> = \>
$[m_{13}-1,m_{23}-2,m_{33}+1,m_{43}+2]$,\\[2mm] 
\>\>$[m]_{14}$ \> 
= \> $[m_{13}-2,m_{23}-1,m_{33}+1,m_{43}+2]$,\\[2mm] 
\>\>$[m]_{15}$ \> = \>
$[m_{13}-2,m_{23}-2,m_{33}+2,m_{43}+2]$.~~~~~~~~~~~~~~~~~(3.26) 
\end{tabbing}
\vspace*{2mm}
In the latest formula (3.26), with the exception of
$[m]_{6}$ and $[m]_{9}$ where a degeneration is present, we 
skip the subscript $k$ in the r.h.s.. The proofs of (3.25) 
and (3.26) follow from direct computations.\\
     
   Using the rule (3.11) which now reads
\begin{eqnarray*}
~~~~~~~~~(m)_{k}&=&\left({[m_{11}-m_{22}]![m_{31}-m_{42}]!
 \over [m_{12}-m_{22}]![m_{12}-m_{11}]!
[m_{32}-m_{42}]![m_{32}-m_{31}]!}\right)^{1/2}\\[5mm] 
&&\times 
(E_{21})^{m_{12}-m_{11}}(E_{43})^{m_{32}-m_{31}}(M)_{k} 
~~~~~~~~~~~~~~~~~~~~~~~~~~~~~~~~~~~~~~~(3.11') 
\end{eqnarray*}\\[2mm]
we can find all the basis vectors $(m)_{k}$ : 
\begin{eqnarray*}
(m)_{0}& = &\left|0,0,0,0;(m)\right>,\\[4mm] 
(m)_{1}
& = &a_{1}q\left\{-\left({[l_{13}-l_{11}][l_{33}-l_{31}+1] 
\over 
[2l+1][2l'+1]}\right)^{1/2}\left|1,0,0,0;(m)^{+11}\right> 
\right.\\
&&-q^{l'-s'}\left(
{[l_{13}-l_{11}][l_{31}-l_{43}-1] \over 
[2l+1][2l'+1]}\right)^{1/2}\left|0,1,0,0;(m)^{+11-31}\right>\\ 
&&+p^{-l+s}\left(
{[l_{11}-l_{23}][l_{33}-l_{31}+1]\over 
[2l+1][2l'+1]}\right)^{1/2}\left|0,0,1,0;(m)\right>\\ 
&   &\left. +p^{-l+s}q^{l'-s'}\left( 
{[l_{11}-l_{23}][l_{31}-l_{43}-1]\over
[2l+1][2l'+1]}\right)^{1/2}\left|0,0,0,1;(m)^{-31}\right> 
\right\},\\[4mm]
 (m)_{2}
& = &-a_{2}q\left ({q\over p}\right)^{l-s-1}
\left\{ \left( {[l_{11}-l_{23}][l_{33}-l_{31}+1]\over 
[2l][2l'+1]}\right)^{1/2}\left|1,0,0,0;(m)^{+11}\right>\right. 
\\
& &+q^{l'-s'}\left(
{[l_{11}-l_{23}][l_{31}-l_{43}-1] \over 
[2l][2l'+1]}\right)^{1/2}\left|0,1,0,0;(m)^{+11-31}\right>\\ 
&  &+ q^{l+s+1}\left(
{[l_{13}-l_{11}][l_{33}-l_{31}+1]\over 
[2l][2l'+1]}\right)^{1/2}\left|0,0,1,0;(m)\right>\\ 
&&\left. +q^{l+s+l'-s'+1}\left(
{[l_{13}-l_{11}][l_{31}-l_{43}-1]\over 
[2l][2l'+1]}\right)^{1/2}\left|0,0,0,1;(m)^{-31}\right> 
\right\},\\[4mm]
(m)_{3}
& = &a_{3}\left\{ -q\left ({q\over p}\right)^{l'-s'} 
\left( {[l_{13}-l_{11}][l_{31}-l_{43}-1]\over
[2l+1][2l']}\right)^{1/2}\left|1,0,0,0;(m)^{+11}\right>\right. 
\\
& &+\left ({q\over p}\right)^{2l'}\left( 
{[l_{13}-l_{11}][l_{33}-l_{31}+1] \over 
[2l+1][2l']}\right)^{1/2}\left|0,1,0,0;(m)^{+11-31}\right>\\ 
&  &+ qp^{-l+s}\left ({q\over p}\right)^{l'-s'}\left( 
{[l_{11}-l_{23}][l_{31}-l_{43}-1]\over 
[2l+1][2l']}\right)^{1/2}\left|0,0,1,0;(m)\right>\\
&&\left. -p^{-l+s}\left ({q\over p}\right)^{2l'} 
\left(
{[l_{11}-l_{23}][l_{33}-l_{31}+1]\over 
[2l+1][2l']}\right)^{1/2}\left|0,0,0,1;(m)^{-31}\right> 
\right\},\\[4mm]
(m)_{4}
& = &a_{4}\left ({q\over p}\right)^{l-s+l'-s'-1} 
\left\{ q\left( {[l_{11}-l_{23}][l_{31}-l_{43}-1]\over
[2l][2l']}\right)^{1/2}\left|1,0,0,0;(m)^{+11}\right>\right. 
\\
& &-p^{-l'-s'}\left(
{[l_{11}-l_{23}][l_{33}-l_{31}+1] \over 
[2l][2l']}\right)^{1/2}\left|0,1,0,0;(m)^{+11-31}\right>\\ 
&  &+ q^{l+s+2}\left(
{[l_{13}-l_{11}][l_{31}-l_{43}-1]\over 
[2l][2l']}\right)^{1/2}\left|0,0,1,0;(m)\right>\\ 
&&\left. -q^{l+s+1}p^{-l'-s'}\left(
{[l_{13}-l_{11}][l_{33}-l_{31}+1]\over 
[2l][2l']}\right)^{1/2}\left|0,0,0,1;(m)^{-31}\right> 
\right\},\\[4mm]
(m)_{5}& = &a_{5}\left ({q\over p}\right)^{l'-s'+1} 
\left\{\left(
{[l_{13}-l_{11}][l_{13}-l_{11}-1] \over 
[2l+1][2l+2]}\right)^{1/2}\left|1,1,0,0;(m)^{+11+11-31} 
\right>\right.
\\ &   &+p^{-l+s}\left( {[l_{13}-l_{11}][l_{11}-l_{23}+1] 
\over
[2l+1][2l+2]}\right)^{1/2}\left|0,1,1,0;(m)^{+11-31}\right>\\ 
& &-p^{-l+s+1}\left( {[l_{13}-l_{11}][l_{11}-l_{23}+1] \over 
[2l+1][2l+2]}\right)^{1/2}\left|1,0,0,1;(m)^{+11-31}\right>\\ 
&&\left. +p^{2(-l+s)}\left( {[l_{11}-l_{23}][l_{11}-l_{23}+1] 
\over
[2l+1][2l+2]}\right)^{1/2}\left|0,0,1,1;(m)^{-31}\right>\right\}, 
\\[4mm]
(m)_{6}
& =& {a_6\over a_5}\left({[2l+1][2l+2][l_{13}-l_{11}] 
\over [l_{11}-l_{23}+1]}\right)^{1/2}(m)_5
-a_{6}\left ({q\over p}\right)^{l'-s'+1}
{[2l+2]\over [2l]}\times \\
& &\times \left\{[l_{13}-l_{11}-2]
\left({[l_{13}-l_{11}-1]\over
[l_{11}-l_{23}+1]}\right)^{1/2}\left|1,1,0,0;(m)^{+11+11-31} 
\right>\right.\\
& &+p^{-l+s+1}[l_{13}-l_{11}-1]
\left|0,1,1,0;(m)^{+11-31}\right>\\
& &-p^{-l+s+2}[l_{13}-l_{11}-1]
\left|1,0,0,1;(m)^{+11-31}\right>\\
&   &\left. +p^{-2(l-s-1)}\left([l_{11}-l_{23}] 
[l_{13}-l_{11}]\right)^{1/2}\left|0,0,1,1;(m)^{-31} 
\right>\right\},\\[4mm]
(m)_{7}
& =& {a_7\over a_5}\left({[2l-1][2l] 
[2l+1][2l+2][l_{13}-l_{11}-1]
[l_{13}-l_{11}]
\over [l_{11}-l_{23}][l_{11}-l_{23}+1]}\right)^{1/2}(m)_5\\ 
&&-{a_7\over a_6}[2][2l+1]\left({[2l-1][l_{13}-l_{11}-1] 
\over [2l][l_{11}-l_{23}]}\right)^{1/2}(m)_6\\ 
&&-a_{7}{[2l+1][2l+2]\over ([2l-1][2l])^{1/2}}
\left ({q\over p}\right)^{l'-s'+1}
\left\{{[l_{13}-l_{11}-2][l_{13}-l_{11}-3]\over 
([l_{11}-l_{23}][l_{11}-l_{23}+1])^{1/2}} 
\left|1,1,0,0;(m)^{+11+11-31}\right>\right.\\
& &+p^{-l+s+2}[l_{13}-l_{11}-2]\left({[l_{13}-l_{11}-1] 
\over [l_{11}-l_{23}]}\right)^{1/2}
\left|0,1,1,0;(m)^{+11-31}\right>\\
& &-p^{-l+s+3}[l_{13}-l_{11}-2]\left({[l_{13}-l_{11}-1] 
\over [l_{11}-l_{23}]}\right)^{1/2}
\left|1,0,0,1;(m)^{+11-31}\right>\\
&   &\left. +p^{-2(l-s-2)}\left([l_{13}-l_{11}-1] 
[l_{13}-l_{11}]\right)^{1/2}\left|0,0,1,1;(m)^{-31} 
\right>\right\},
\end{eqnarray*}
\begin{eqnarray*}
(m)_{8}
& = &a_{8}q^2\left({q\over p}\right)^{l-s-1}\left\{\left( 
{[l_{33}-l_{31}+1][l_{33}-l_{31}+2] \over 
[2l'+1][2l'+2]}\right)^{1/2}\left|1,0,1,0;(m)^{+11}\right> 
\right.\\
       &   &+q^{2l'}\left(
{[l_{33}-l_{31}+2][l_{31}-l_{43}-1] \over 
[2l'+1][2l'+2]}\right)^{1/2}\left|1,0,0,1;(m)^{+11-31}\right>\\ 
& &+q^{2l'+1}\left( {[l_{33}-l_{31}+2][l_{31}-l_{43}-1]
\over
[2l'+1][2l'+2]}\right)^{1/2}\left|0,1,1,0;(m)^{+11-31}\right>\\ 
& &\left. +q^{2(2l'+1)}\left( {[l_{31}-l_{43}-2] 
[l_{31}-l_{43}-1]\over
[2l'+1][2l'+2]}\right)^{1/2}\left|0,1,0,1;(m)^{+11-31-31} 
\right>\right\},\\[4mm]
(m)_{9}
& = &{a_9\over a_8}
\left({[2l'+1][2l'+2][l_{33}-l_{31}+2]\over 
[l_{31}-l_{43}-1]}\right)^{1/2}(m)_8\\
&&-a_{9}q^2\left({q\over p}\right)^{l-s-1} 
{[2l'+2]\over [2l']}
\left\{[l_{33}-l_{31}]\left(
{[l_{33}-l_{31}+1]\over
[l_{31}-l_{43}-1]}\right)^{1/2}\left|1,0,1,0;(m)^{+11}\right> 
\right.\\
&   &+q^{l'-s'-1}[l_{33}-l_{31}+1] 
\left|1,0,0,1;(m)^{+11-31}\right>\\
& &+q^{l'-s'}[l_{33}-l_{31}+1] 
\left|0,1,1,0;(m)^{+11-31}\right>\\
& &\left. +q^{2(l'-s')}
\left([l_{33}-l_{31}+2][l_{31}-l_{43}-2]\right)^{1/2} 
\left|0,1,0,1;(m)^{+11-31-31}
\right>\right\},\\[4mm]
(m)_{10}
& = &{a_{10}\over a_8}\left({[2l'-1][2l'][2l'+1][2l'+2] 
[l_{33}-l_{31}+1][l_{33}-l_{31}+2]\over 
[l_{31}-l_{43}-2][l_{31}-l_{43}-1]}\right)^{1/2}(m)_8\\
&&-{a_{10}\over a_9}[2][2l'+1]
\left({[2l'-1][l_{33}-l_{31}+1]\over 
[2l'][l_{31}-l_{43}-2]}\right)^{1/2}(m)_9\\ 
&&-a_{10}\left({q\over p}\right)^{l-s-1} 
{
[2l'+1][2l'+2]\over [2l'-1][2l']} 
\left\{q^2[l_{33}-l_{31}-1][l_{33}-l_{31}]\times \right.\\ 
&&\times \left({[2l'-1][2l']\over
[l_{31}-l_{43}-2][l_{31}-l_{43}-1]}\right)^{1/2} 
\left|1,0,1,0;(m)^{+11}\right>\\
&   &+q^{l'-s'}[l_{33}-l_{31}] 
\left({[2l'-1][2l'][l_{33}-l_{31}+1]\over 
[l_{31}-l_{43}-2]}\right)^{1/2}
\left|1,0,0,1;(m)^{+11-31}\right>\\
& &+q^{l'-s'+1}[l_{33}-l_{31}] 
\left({[2l'-1][2l'][l_{33}-l_{31}+1]\over 
[l_{31}-l_{43}-2]}\right)^{1/2}
\left|0,1,1,0;(m)^{+11-31}\right>\\
& &\left. +q^{2(l'-s')}
\left([2l'-1][2l'][l_{33}-l_{31}+1]
[l_{33}-l_{31}+2]\right)^{1/2}
\left|0,1,0,1;(m)^{+11-31-31}
\right>\right\},\\[4mm]
(m)_{11} & = &a_{11}\left({q\over p}\right)^{l-s+l'-s'+1} 
\left\{p\left(
{[l_{13}-l_{11}-1][l_{33}-l_{31}+2] \over 
[2l+1][2l'+1]}\right)^{1/2}\left|1,1,1,0;(m)^{+11+11-31} 
\right>\right.\\
&  &+p^{-l+s+2}\left( {[l_{11}-l_{23}+1][l_{33}-l_{31}+2] 
\over
[2l+1][2l'+1]}\right)^{1/2}\left|1,0,1,1;(m)^{+11-31} 
\right>\\
& &+q^{l'-s'+1}p\left( {[l_{13}-l_{11}-1][l_{31}-l_{43}-2] 
\over [2l+1][2l'+1]}\right)^{1/2}\left|1,1,0,1;(m)^{+11+11-31-31} 
\right>\\
& &\left. +q^{l'-s'+2}p^{-l+s+1}\left( 
{[l_{11}-l_{23}+1][l_{31}-l_{43}-2] \over 
[2l+1][2l'+1]}\right)^{1/2}\left|0,1,1,1;(m)^{+11-31-31} 
\right>\right\},\\[4mm]
(m)_{12} & = &{a_{12}\over a_{11}}
\left({[2l][2l+1][l_{13}-l_{11}-1]\over 
[l_{11}-l_{23}+1]}\right)^{1/2} (m)_{11}\\ 
&&-a_{12}p\left({q\over p}\right)^{l-s+l'-s'} 
{[2l+1]\over [2l]}\left\{[l_{13}-l_{11}-2] 
\left({[2l][l_{33}-l_{31}+2] \over [l_{11}-l_{23}+1] 
[2l'+1]}\right)^{1/2}\times \right.\\
&&\times
\left|1,1,1,0;(m)^{+11+11-31}
\right>\\
&  &+p^{-l+s+2}\left( {[2l][l_{13}-l_{11}-1][l_{33}-l_{31}+2] 
\over
[2l'+1]}\right)^{1/2}\left|1,0,1,1;(m)^{+11-31} 
\right>\\
& &+q^{l'-s'+1}[l_{13}-l_{11}-2] 
\left({[2l][l_{31}-l_{43}-2]
\over [[l_{11}-l_{23}+1]][2l'+1]}\right)^{1/2} 
\left|1,1,0,1;(m)^{+11+11-31-31}
\right>\\
& &\left. +q^{l'-s'+2}p^{-l+s+1}\left( 
{[2l][l_{13}-l_{11}-1][l_{31}-l_{43}-2] \over 
[2l'+1]}\right)^{1/2}\left|0,1,1,1;(m)^{+11-31-31} 
\right>\right\},\\[4mm]
(m)_{13} & = &{a_{13}\over a_{11}}
\left({[2l'][2l'+1][l_{33}-l_{31}+2]\over 
[l_{31}-l_{43}-2]}\right)^{1/2} (m)_{11}\\ 
&&-a_{13}p\left({q\over p}\right)^{l-s+l'-s'} 
{[2l'+1]\over [2l']}\left\{[l_{33}-l_{31}+1] 
\left({[2l'][l_{13}-l_{11}-1] \over [2l+1] 
[l_{31}-l_{43}-2]}\right)^{1/2}\times \right.\\ 
&&\times
\left|1,1,1,0;(m)^{+11+11-31}
\right>\\
&  &+p^{-l+s+1}[l_{33}-l_{31}+1]
\left( {[l_{11}-l_{23}+1][2l']
\over
[2l+1][l_{31}-l_{43}-2]}\right)^{1/2}\left|1,0,1,1;(m)^{+11-31} 
\right>\\
& &+q^{l'-s'}
\left({[l_{13}-l_{11}-1][l_{33}-l_{31}+2][2l'] 
\over [2l+1]}\right)^{1/2}
\left|1,1,0,1;(m)^{+11+11-31-31}
\right>\\
& &\left. +q^{l'-s'+1}p^{-l+s}\left( 
{[l_{11}-l_{23}+1][l_{33}-l_{31}+2][2l']\over 
[2l+1]}\right)^{1/2}\left|0,1,1,1;(m)^{+11-31-31} 
\right>\right\},\\[4mm]
(m)_{14} & = &{a_{14}\over a_{11}}
\left({[2l][2l+1][l_{13}-l_{11}-1][2l'][2l'+1] 
[l_{33}-l_{31}+2]\over [l_{11}-l_{23}+1] 
[l_{31}-l_{43}-2]}\right)^{1/2} (m)_{11}\\ 
&&-{a_{14}\over a_{12}}
\left({[2l'][2l'+1][l_{33}-l_{31}+2]\over 
[l_{31}-l_{43}-2]}\right)^{1/2} (m)_{12}\\ 
&&-{a_{14}\over a_{13}}
\left({[2l][2l+1][l_{13}-l_{11}-1]\over 
[l_{11}-l_{23}+1]}\right)^{1/2} (m)_{13}\\ 
&&-a_{14}p\left({q\over p}\right)^{l-s+l'-s'-1} 
{[2l+1][2l'+1]\over [2l][2l']}
\left\{ 
[l_{13}-l_{11}-2][l_{33}-l_{31}+1]\times \right.\\ 
&&\times 
\left({[2l][2l']\over [l_{11}-l_{23}+1] 
[l_{31}-l_{43}-2]}\right)^{1/2}
\left|1,1,1,0;(m)^{+11+11-31}
\right>\\
&  &+p^{-l+s+2}[l_{33}-l_{31}+1]
\left( {[l_{13}-l_{11}-1][2l][2l']
\over
[l_{31}-l_{43}-2]}\right)^{1/2}\left|1,0,1,1;(m)^{+11-31} 
\right>\\
& &+q^{l'-s'}[l_{13}-l_{11}-2] 
\left({[2l][2l'][l_{33}-l_{31}+2]
\over [l_{11}-l_{23}+1]}\right)^{1/2} 
\left|1,1,0,1;(m)^{+11+11-31-31}
\right>\\
& &\left. +q^{l'-s'+1}p^{-l+s+1}
\left([2l][l_{13}-l_{11}-1][l_{33}-l_{31}+2][2l'] 
\right)^{1/2}\left|0,1,1,1;(m)^{+11-31-31} 
\right>\right\},\\[4mm]
(m)_{15}&=& a_{15}(m)
~~~~~~~~~~~~~~~~~~~~~~~~~~~~~~~~~~~~~~~~~~ 
~~~~~~~~~~~~~~~~~~~~~~~~~~~~~~~~~~~~(3.27) 
\end{eqnarray*}
\\
where $l={1\over 2}(m_{13}-m_{23})$, $s=m_{11}-{1\over 
2}(m_{13}+m_{23})$, $l'={1\over 2}(m_{33}-m_{43})$  and 
$s'=m_{31}-{1\over 2}(m_{33}+m_{43})$, while $(m)_k^{\pm ij}$ 
is a GZ basis vector obtained from $(m)_k$ with replacing the 
element
$m_{ij}$ by $m_{ij}\pm 1$. We can write down the coefficients 
in (3.27) all in terms of $l$, $s$, $l'$ and $s'$ only but 
here we leave them partially expressed in terms of $l_{ij}$ 
and $l'_{ij}$. From (3.27) we can immediately find
all the (local) lowest weight vectors $(M)^{V}_{k}$ which, by 
definition, are annihilated by $E_{21}$ and $E_{43}$. Let us 
remind again that every firdirmod $V_k^{p,q}$ on a level $\eta$, 
spanned on linear combinations of
$\left|\theta_1,\theta_2,\theta_3,\theta_4; (m)\right>$ 
in (3.14) with a fixed
$\sum_{i=1}^4\theta_i\equiv \eta$ is a linear space of 
homogeneous
polynomials of a definite power $\eta$ in the
negative odd generators $E_{ij}$ ($1\leq j< 3\leq i\leq 4$) 
acting on $(m)\in V_0^{p,q}([m])$.
Taking into account all results obtained above we have 
proved the following assertion\\[5mm]
{\bf Proposition} 3 : {\it Every $U_{p,q}[gl(2/2)_{0}]$--fidirmod 
$V^{p,q}_{k}$ in decomposition (3.16) is characterized by
a highest weight $[m]_{k}$ given in (3.25) and is spanned 
by a GZ basis $(m)_{k}$ given in (3.27)}.\\
         
  The latest formula (3.27), in fact, represents a way in which
the reduced basis is expressed in terms of the induced basis and 
vas versa it is not a problem for us to find the inverse relation 
between these bases (see the Appendix). For further convenience 
the vectors $(m)_{\tilde k}\equiv (m)_k$ (for $k=6,7,9,10,12,13$ 
and 14) are partially given via other $(m)_k$ which are 
completely expressed in terms of
$\left|\theta_1,\theta_2,\theta_3,\theta_4; (m)\right>$. 
It is not difficult to write down the explicit
decompositions of these $(m)_{\tilde k}$ in the induced basis. 
But here we prefer the expressions in (3.27) which are more 
compact and more convenient for finding the inverse relation 
between two bases and matrix elements of odd generators.\\
     
  Now we are ready to canculate the matrix elements of the
generators $E_{ij}$. It is sufficient to canculate the 
matrix elements of the Cartan--Chevalley generators only, 
since $U_{p,q}[gl(2/2)]$ can be generated by these 
generators and any its representation in some basis is 
completely defined by their actions on the same basis. 
For the even generators which do not shift the $\eta$'s we 
readily have
\begin{eqnarray*}
~~~~~~~~~~~~~~E_{11}(m)_{k}& = &(l_{11}+1)(m)_{k},\\ 
E_{22}(m)_{k}& = &(l_{12}+l_{22}-l_{11}+2)(m)_{k},\\ 
E_{12}(m)_{k}& =
&\left([l_{12}-l_{11}][l_{11}-l_{22}]\right)^{1/2}(m)_{k}^{+11},\\ 
E_{21}(m)_{k}& =
&\left([l_{12}-l_{11}+1][l_{11}-l_{22}-1] 
\right)^{1/2}(m)_{k}^{-11},\\
L_1(m)_k&=&{1\over 2}(l_{12}-l_{22}-1)(m)_k,\\ 
E_{33}(m)_{k}& = & (l_{31}+1)(m)_{k},\\ 
E_{44}(m)_{k}& = 
&(l_{32}+l_{42}-l_{31}+2)(m)_{k},\\ 
E_{34}(m)_{k}& = &
\left([l_{32}-l_{31}][l_{31}-l_{42}]\right)^{1/2}(m)_{k}^{+31},\\ 
E_{43}(m)_{k}& =&
\left([l_{32}-l_{31}+1][l_{31}-l_{42}-1]\right)^{1/2}
(m)_{k}^{-31},\\
L_3(m)_k&=&{1\over 2}(l_{32}-l_{42}-1)(m)_k.
~~~~~~~~~~~~~~~~~~~~~~~~~~~~~~~~~~~~~~~~~(3.28) 
\end{eqnarray*}
\vspace*{2mm}
As the matrix elements of $E_{23}$ and $E_{32}$ 
are very long expressions 
we only explain here how to find them. 
By construction a reduced basis 
vector $(m)_k$ in (3.27) belonging to a fidirmod $V_k^{p,q}$
on a level $\eta$ is a homogeneous polynomial of a power $\eta$ 
in odd generators
$E_{ij}, ~1\leq j< 3\leq i\leq 4$ acting on $(m)\in V_0^{p,q}$. 
Under the action of  $E_{23}$ (or $E_{32}$, respectively) this 
vector $(m)_k$ is shifted to other fidirmods $V^{p,q}_{k'}$ on 
the previous level $\eta-1$ (or on the next level $\eta +1$, 
respectively), i.e., we get on the r.h.s of $E_{23}(m)_k$
(or $E_{32}(m)_k$, respectively) a homogeneous polynomial
of a degree $\eta-1$ (or $\eta+1$, respectively). Using the 
inverse relations (A.1) we can express the latter polynomials 
obtained in terms of the reduced basis, that is we get 
matrix elements of $E_{23}$ and $E_{32}$ in this basis.
It is a standard way to find matrix elements but in 
practice we can use a trick making calculations simpler. 
Since $E_{23}$ commutes with $E_{21}$ and $E_{43}$ we 
first calculate the action of $E_{23}$ on the highest
vectors only and then apply (3.11) to find all matrix 
elements of this generator on arbitrary $(m)_k$. It is 
less complicated 
to compute matrix elements of $E_{32}$ in the standard way 
but we can apply a similar trick, namely, we first calculate 
the action of $E_{23}$ (which commutes with $E_{12}$ and 
$E_{34}$) on the lowest weight vectors and then apply the 
rule inverse to (3.11).\\
     
  It can be shown that the representations constructed contain
all finite--dimensional irreducible representations of 
$U_{p,q}[gl(2/2)]$ classified as typical or nontypical 
representations which are subjects of next investigations.\\[5mm] 
{\large {\bf IV. Conclusion}}
     
  We have considered the two--parametric quantum 
deformations $U_{p,q}[gl(2/2)]$ and described in 
detail a method for constructing its 
finite--dimensional representations. The 
representations constructed can be decomposed into 
finite--dimensional irreducible representations of 
the even subalgebra  $U_{p,q}[gl(2/2)_0]$ and 
therefore can be given in bases of the latter. Using 
Poincar\'e-Birkhoff-Witt theorem and the induced 
representation method we constructed the induced 
basis of the induced module $W^{p,q}$. This basis, 
however, does not allow a clear description of a 
decomposition of $W^{p,q}$ into 
$U_{p,q}[gl(2/2)_0]$--fidirmods. It was the reason 
the reduced basis was introduced. The latter basis 
representing a union of GZ bases of the even 
subalgebra $U_{p,q}[gl(2/2)_0]$ according to the 
branching rule $U_{p,q}[gl(2/2)]\supset 
U_{p,q}[gl(2/2)_{0}] \supset gl(1)\otimes gl(1)$
is refered to as Quasi--GZ basis. This step is 
intermediate but of independent interest. Having 
these two bases, the induced and the reduced ones, 
and the relations between them we can find all matrix
elements of finite--dimensional representations of 
$U_{p,q}[gl(2/2)]$. It turn out that the representations 
constructed contain all finite--dimensional irreducible 
representations of $U_{p,q}[gl(2/2)]$ and can be 
classified into typical and nontypical representations 
which are subjects of later papers.\\
     
  Looking at the basis transformations and the matrix 
elements we observe, even at generic deformation 
parameters, some "anomalies" which are canceled out 
at $p=q$. It means that the finite--dimensional 
representations of the two-parametric quantum 
superalgebra $U_{q}[gl(2/2)]$ are not simply trivial 
deformations from those of the classical Lie superalgebra 
$gl(2/2)$ in the sense that they can not be found from 
classical analogues by putting quantum deformation brackets 
in appropriate places unlike many cases of one--parametric 
deformations. For example, the expressions
$$ {q\over p}[2l][l_{13}-l_{11}] - 
[2l+1][l_{13}-l_{11}-1] \eqno(5.1)$$
and
$${q\over p} [2l'][l_{33}-l_{31}] - 
[2l'+1][l_{33}-l_{31}-1] \eqno(5.2)$$
appearing in the basis transformations and matrix 
elements can be written in the forms:
$$({q\over p}-1)[2l][l_{13}-l_{11}]+
\left({q\over p}\right)
^{l_{13}-l_{11}-1}[l_{11}-l_{23}]\eqno(5.1')$$ 
and
$$({q\over p}-1)[2l'][l_{33}-l_{31}]+\left({q\over p}\right) 
^{l_{33}-l_{31}-1}[l_{31}-l_{43}],\eqno(5.2')$$ 
respectively. At $p=q$ the latest expressions become 
$[l_{11}-l_{23}]$ and $[l_{31}-l_{43}]$, respectively, 
exactly as in the one--parametric case \cite{ky2,ky3}.\\

  We hope that it is not very difficult to extend the
present method to the case of one or both deformation 
parameters being roots of unity. For conclusion, let us 
emphasize that our method has an advantage that it avoids 
the use of the Clebsch--Gordan coefficients which are not 
always known, especially for higher rank (classical and 
quantum) groups and multi--parametric deformations.
\newpage
{\bf Acknowledgements}\\

 I would like to thank the Nishina memorial foundation 
for financial support and the Department of Physics, Chuo 
University, Tokyo, Japan for warm hospitality. Fruitful 
discussions with K. Furuta, T. Inami and other members 
of the Theory Group of the Department of Physics, Chuo 
University, are also hereby acknowledged.\\     
     
  This work was supported in part by the National Research
Programme for Natural Sciences of Vietnam under grant  
number KT -- 04.1.1.

\newpage
\begin{flushleft}
{\bf Appendix}\\
\end{flushleft}
     
  The induced basis (4.20) is expressed in terms 
of the reduced basis through the following inverse 
relation 
\begin{eqnarray*}
\left|1,0,0,0;(m)\right>& = &
-{1\over a_{1}}q^{l+s-1}
p^{-l'-s'}
\left({[l_{13}-l_{11}+1][l_{33}-l_{31}+1]\over 
[2l+1][2l'+1]} \right)^{1/2}(m)_{1}^{-11}\\[2mm]
&   &-{1\over
a_{2}}{q^{-l+s-1}p^{-l'-s'-1}\over [2l+1]} 
\left({[2l][l_{11}-l_{23}-1][l_{33}-l_{31}+1]\over [2l'+1]} 
\right)^{1/2}(m)_{2}^{-11}\\[2mm]
&   &-{1\over
a_{3}}{q^{l+s}p^{l'-s'}\over [2l'+1]} 
\left({[l_{13}-l_{11}+1][l_{31}-l_{43}-1][2l']\over [2l+1]} 
\right)^{1/2}(m)_{3}^{-11}\\[2mm]
&   & +{1\over
a_{4}}{q^{-l+s}p^{l'-s'-1}\over [2l+1][2l'+1]}
\left([2l][l_{11}-l_{23}-1][2l'][l_{31}-l_{43}-1] 
\right)^{1/2}(m)_{4}^{-11},\\[2mm]
\left|0,1,0,0;(m)\right>& = &-{1\over a_{1}}
q^{l+s}
\left({[l_{13}-l_{11}+1][l_{31}-l_{43}]\over [2l+1][2l'+1]} 
\right)^{1/2}(m)_{1}^{-11+31}\\[2mm]
&   &-{1\over
a_{2}}{q^{-l+s}\over p[2l+1]}
\left({[2l][l_{11}-l_{23}-1][l_{31}-l_{43}]\over [2l'+1]} 
\right)^{1/2}(m)_{2}^{-11+31}\\[2mm]
&   &+{1\over
a_{3}}\left({p\over q}\right)^{l'-s'-1} 
{q^{l+s}\over [2l'+1]} 
\left({[l_{13}-l_{11}+1][2l'][l_{33}-l_{31}]\over 
[2l+1]} \right)^{1/2}(m)_{3}^{-11+31}\\[2mm]
&   &
-{1\over a_{4}}\left({p\over q}\right)^{l'-s'-2}
{q^{-l+s-1}\over [2l+1][2l'+1]}
\left({[2l][l_{11}-l_{23}-1][2l']\over [l_{33}-l_{31}} 
\right)^{1/2}
(m)_{4}^{-11+31},\\[2mm]
\left|0,0,1,0;(m)\right>& = &-{1\over a_{1}}q^{-1} 
p^{-l'-s'}\left({[l_{11}-l_{23}][l_{33}-l_{31}+1]
\over [2l+1][2l'+1]}\right)^{1/2}(m)_{1}\\[2mm]
&   &+{1\over
a_{2}}\left({p\over q}\right)^{l-s}
{p^{-l'-s'-1}\over [2l+1]} 
\left({[2l][l_{13}-l_{11}][l_{33}-l_{31}+1]\over [2l'+1]} 
\right)^{1/2}(m)_{2}\\[2mm]
&   &+{1\over a_{3}}{p^{l'-s'}\over [2l'+1]} 
\left({[l_{11}-l_{23}][2l'][l_{31}-l_{43}-1]\over [2l+1]} 
\right)^{1/2}(m)_{3}\\[2mm]
&   &+{1\over
a_{4}}\left({p\over q}\right)^{l-s-1} 
{p^{l'-s'}\over [2l+1][2l'+1]}\times \\
& &\times  
\left([2l][l_{13}-l_{11}][2l'][l_{31}-l_{43}-1] 
\right)^{1/2}(m)_{4},\\[2mm]
\left|0,0,0,1;(m)\right>& = &{1\over a_{1}}
\left({[l_{11}-l_{23}][l_{31}-l_{43}+1]\over [2l+1][2l'+1]} 
\right)^{1/2}(m)_{1}^{+31}\\[2mm]
&   &-{1\over a_{2}}\left({p\over q}\right)^{l-s-1} 
{1\over [2l+1]}
\left({[2l][l_{13}-l_{11}][l_{31}-l_{43}] 
\over [2l'+1]}
\right)^{1/2}(m)_{2}^{+31}\\[2mm]
&   &-{1\over a_{3}}\left({p\over q}\right)^{l'-s'-1} 
{1\over [2l'+1]}
\left({[l_{11}-l_{23}][2l'][l_{33}-l_{31}]\over [2l+1]} 
\right)^{1/2}(m)_{3}^{+31}\\[2mm]
&   & -{1\over a_{4}}\left({p\over q}\right)^{l-s+l'-s'-2} 
{1\over [2l+1][2l'+1]}\times \\
& &\times  
\left({[2l][2l']\over [l_{13}-l_{11}][l_{31}-l_{43}]} 
\right)^{1/2}(m)_{4}^{+31},\\[2mm]
\left|1,1,0,0;(m)\right>& = &{1\over a_{5}} 
q^{2(l+s)}\left({p\over q}\right)^{l'-s'} 
\left({[l_{13}-l_{11}+1][l_{13}-l_{11}+2]\over [2l+1][2l+2]} 
\right)^{1/2}(m)_{5}^{-11-11+31}\\[2mm]
&   &+{1\over a_{6}}
\left({p\over q}\right)^{l'-s'-3}
{q^{-2(l-s)-3}\over [2l+2]} 
\left( p^2q^{2l-1}+(p-q)[2l-1]\right)\times \\
&&\times \left([l_{11}-l_{23}-1][l_{13}-l_{11}-1]\right)^{1/2} 
(m)_{6}^{-11-11+31}\\[2mm]
&   &-{1\over
a_{7}} 
\left({p\over q}\right)^{l'-s'-3}
{q^{-2(l-s+1)}\over [2][2l+1][2l+2]}
\times \\ 
&&\times
\left([2l][2l-1][l_{11}-l_{23}-2][l_{11}-l_{23}-1] 
\right)^{1/2}(m)_{7}^{-11-11+31},\\[2mm]
E^{(-2)}&\equiv &
\left|0,1,1,0;(m)\right> - p\left|1,0,0,1;(m)\right>\\ 
& = &{1\over a_{5}}q^{l+s+1}
\left({p\over q}\right)^{l'-s'+1}[2] 
\left({[l_{13}-l_{11}+1][l_{11}-l_{23}]\over [2l+1][2l+2]} 
\right)^{1/2}(m)_{5}^{-11+31}\\[2mm]
&   &+{1\over a_{6}}\left({p\over q}\right)^{l-s+l'-s'-2}
{q^{-l+s}\over [2l+2]} 
\left\{[2][2l-1][l_{13}-l_{11}] \right.\\
&&-[2l]\left(p^2[l_{13}-l_{11}] 
+q^{-1}[l_{13}-l_{11}-1] \right)\left.
\right\}(m)_{6}^{-11+31}\\[2mm] 
&&+{1\over a_{7}}
\left({p\over q}\right)^{l-s+l'-s'-2}
{q^{-l+s}\over [2l+1][2l+2]}
\times \\ 
&&\times
\left([2l][2l-1][l_{13}-l_{11}][l_{11}-l_{23}-1] 
\right)^{1/2}(m)_{7}^{-11+31},
\end{eqnarray*}
\begin{eqnarray*}
\left|0,0,1,1;(m)\right>& = &{1\over a_{5}} 
\left({p\over q}\right)^{l'-s'}
\left({[l_{11}-l_{23}][l_{11}-l_{23}+1]\over [2l+1][2l+2]} 
\right)^{1/2}(m)_{5}^{+31}\\[2mm]
&   &+{1\over a_{6}}\left({p\over q}\right)^{2(l-s)+l'-s'-4}
\left\{ {p\over q}[2l][l_{13}-l_{11}-2]
-[2l-1][l_{13}-l_{11}-1]\right\}\times\\
&&\times 
{1\over [2l+2]}\left({[l_{13}-l_{11}]\over [l_{11}-l_{23}]} 
\right)^{1/2}(m)_{6}^{+31}\\[2mm]
&   &-{1\over a_{7}}
\left({p\over q}\right)^{2(l-s)+l'-s'-4} 
{\left([2l][2l-1][l_{13}-l_{11}-1][l_{13}-l_{11}]\right)^{1/2} 
\over [2][2l+1][2l+2]}(m)_{7}^{+31},\\[2mm]
\left|1,0,1,0;(m)\right>& = &{1\over a_{8}}q^{-1} 
p^{-2(l'+s')}\left({q\over p}\right)^{-l+s} 
\left({[l_{33}-l_{31}+1][l_{33}-l_{31}+2]\over [2l'+1][2l'+2]} 
\right)^{1/2}(m)_{8}^{-11}\\[2mm]
&   &+{1\over a_9}q^{-1}p^{2(l'-s')-1}\left({p\over q}\right)^{l-s} 
{\left([l_{33}-l_{31}+1][l_{31}-l_{43}-1]\right)^{1/2}
\over [2l'+2]}(m)_{9}^{-11}\\[2mm]
&   &
-{1\over a_{10}}\left({p\over q}\right)^{l-s}
{qp^{2(l'-s')-1}\over [2][2l'+1][2l'+2]}\times \\
& &\times  
\left([2l'][2l'-1][l_{31}-l_{43}-2][l_{31}-l_{43}-1]\right)^{1/2} 
(m)_{10}^{-11},\\[2mm]
E^{(+2)}&\equiv &
\left|0,1,1,0;(m)\right> +q^{-1}\left|1,0,0,1;(m)\right>\\ 
& = &
{1\over a_{8}}
q^{-1}p^{-l'-s'-1}
\left({p\over q}\right)^{l-s+1}[2]
\left({[l_{33}-l_{31}+1][l_{31}-l_{43}]\over 
[2l'+1][2l'+2]} \right)^{1/2}(m)_{8}^{-11+31}\\[2mm] 
& &
+{1\over a_{9}}
\left({p\over q}\right)^{l-s+l'-s'-2}
{p^{l'-s'}\over q[2l'+2]}\times \\ 
&& \times \left\{[2][2l-1][l_{33}-l_{31}]- 
[2l']
\left (q^{-2}[l_{33}-l_{31}] + p[l_{33}-l_{31}-1]\right)
\right\}(m)_{9}^{-11+31}\\[2mm]
&   &
+{1\over a_{10}}
\left({p\over q}\right)^{l-s+l'-s'-2}
{p^{l'-s'}\over 
q[2l'+1][2l'+2]}\times \\
& &\times  
\left([2l][2l-1][l_{33}-l_{31}][l_{31}-l_{43}-1] 
\right)^{1/2}(m)_{10}^{-11+31},\\[2mm] 
\left|0,1,0,1;(m)\right>& = &{1\over a_{8}} 
\left({p\over q}\right)^{l-s}
\left({[l_{31}-l_{43}][l_{31}-l_{43}+1]\over [2l'+1][2l'+2]} 
\right)^{1/2}(m)_{8}^{-11+31+31}\\[2mm]
&   &+{1\over a_{9}}
\left({p\over q}\right)^{l-s+2(l'-s')-4} 
{1\over [2l'+2]}\left({[l_{33}-l_{31}]\over
[2l'+1]} \right)^{1/2}\times \\
& & \times \left({p\over q}[2l'][l_{33}-l_{31}-2]
-[2l'-1][l_{33}-l_{31}-1]\right) (m)_{9}^{-11+31+31} 
\\[2mm]
&   &-{1\over a_{10}}
\left({p\over q}\right)^{l-s+2(l'-s')-4}
{1\over [2][2l'+1][2l'+2]}\times \\
& &\times \left([2l'][2l'-1][l_{33}-
l_{31}-1][l_{33}-l_{31}]\right)^{1/2} 
(m)_{10}^{-11+31+31},\\[2mm]
\left|1,1,1,0;(m)\right>
& = & q^{-l+s-2}p^{l'-s'-2} 
\left({p\over q}\right)^{l-s+l'-s'}\times \\
& &\times 
\left\{{1\over a_{11}}
\left({[l_{13}-l_{11}+1][l_{33}-l_{31}+1]\over [2l+1][2l'+1]} 
\right)^{1/2}(m)_{11}^{-11-11+31}\right.\\[2mm]
&   &+{1\over a_{12}}
{q\over [2l+1]}
\left({[2l][l_{11}-l_{23}-1][l_{33}-l_{31}+1]\over [2l'+1]} 
\right)^{1/2}(m)_{12}^{-11-11+31}\\[2mm]
&   &+{1\over a_{13}}
{q^2\over [2l'+1]}
\left({[l_{13}-l_{11}+1][2l'][l_{31}-l_{43}-1]\over [2l+1]} 
\right)^{1/2}(m)_{13}^{-11-11+31}\\[2mm]
&   &
+\left.{1\over
a_{14}}{q^3\over [2l+1][2l'+1]}
\left([2l][l_{11}-l_{23}-1][2l'][l_{31}-l_{43}-1] 
\right)^{1/2}(m)_{14}^{-11-11+31}\right\},\\[2mm] 
\left|1,0,1,1;(m)\right>& = & 
p^{l'-s'}
\left({p\over q}\right)^{2(l-s)+l'-s'-3} 
\left\{{1\over a_{11}}
q^{-3}\left(q[2l][l_{13}-l_{11}]- p[2l+1][l_{13}-l_{11}-1]\right) 
\times \right.\\
&&\times
\left( p[2l'+1]-q^2[2l'] \right)\left({[l_{33}-l_{31}+1]
\over [2l+1][2l'+1][l_{11}-l_{23}]}\right)^{1/2}
(m)_{11}^{-11+31}\\[2mm]
&   &-{1\over
a_{12}}{q^{-2}\over [2l+1]}
\left( p[2l'+1]-q^2[2l'] \right)\times \\
& &\times 
\left({[2l][l_{13}-l_{11}][l_{33}-l_{31}+1]\over [2l'+1]} 
\right)^{1/2}(m)_{12}^{-11+31}\\[2mm]
&   &-{1\over a_{13}}
{q^{-1}\over [2l'+1]}
\left(q^{-1}[2l][l_{13}-l_{11}]- p[2l+1][l_{13}-l_{11}-1]\right)
\times\\ 
&&\times
\left({[2l'][l_{31}-l_{43}-1]\over
[2l+1][l_{11}-l_{23}]} \right)^{1/2}(m)_{13}^{-11+31}\\[2mm] 
&   &\left.
+{1\over a_{14}}{1\over [2l+1][2l'+1]} 
\left([2l][l_{13}-l_{11}][2l'][l_{31}-l_{43}-1]\right)^{1/2} 
(m)_{14}^{-11+31}\right\},\\[2mm]
\left|1,1,0,1;(m)\right>& = & 
q^{-l+s-3}\left({p\over q}\right)^{l-s+2(l'-s')-4} 
\left\{{1\over a_{11}}
\left(q[2l]-p^2[2l+1]\right)\times \right.\\
&&\times \left(p[2l'+1][l_{33}-l_{31}-1]-q[2l'][l_{33}-l_{31}] 
\right)\\
&&\times
\left({[l_{13}-l_{11}+1]\over [2l+1][2l'+1][l_{31}-l_{43}]} 
\right)^{1/2}(m)_{11}^{-11-11+31+31}\\[2mm]
&   &-{1\over a_{12}}
{q\over [2l+1]}
\left(p[2l'+1][l_{33}-l_{31}-1]-q[2l'][l_{33}-l_{31}] 
\times \right)\\ 
&&\times \left({[2l][l_{11}-l_{23}-1]\over
[2l'+1][l_{31}-l_{43}]} \right)^{1/2}(m)_{12}^{-11-11+31+31}\\[2mm] 
&&+{1\over a_{13}}
{q\over [2l'+1]}
\left(q[2l]-p^2[2l+1]\right)\times \\
& &\times 
\left({[l_{13}-l_{11}+1][2l'][l_{33}-l_{31}]\over [2l+1]} 
\right)^{1/2}(m)_{13}^{-11-11+31+31}\\[2mm]
&   &\left. +{1\over a_{14}}
{q^2\over[2l+1][2l'+1]}
\left([2l][l_{11}-l_{23}-1][2l'][l_{33}-l_{31}] 
\right)^{1/2}(m)_{14}^{-11-11+31+31}\right\},\\[2mm] 
\left|0,1,1,1;(m)\right>& = & 
\left({p\over q}\right)^{2(l-s+l'-s'-2)} 
\left\{{1\over a_{11}}
q^{-2}\left(p[2l+1][l_{13}-l_{11}-1]-q[2l][l_{13}-l_{11}] 
\right)\times\right.\\
&&\times
\left(p[2l'+1][l_{33}-l_{31}-1]-q[2l'][l_{33}-l_{31}] 
\right)\times\\
&&\times
\left({1\over [2l+1][2l'+1][l_{11}-l_{23}][l_{31}-l_{43}]} 
\right)^{1/2}(m)_{11}^{-11+31+31}\\[2mm]
&   &+{1\over a_{12}}
{q^{-1}\over [2l+1]}
\left(p[2l'+1][l_{33}-l_{31}-1]-q[2l'][l_{33}-l_{31}] 
\right)\times\\
&&\times
\left({[2l][l_{13}-l_{11}]\over [2l'+1][l_{31}-l_{43}]} 
\right)^{1/2}(m)_{12}^{-11+31+31}\\[2mm]
&   &+{1\over a_{13}}
{q^{-1}\over [2l'+1]}
\left(p[2l+1][l_{13}-l_{11}-1]-q[2l][l_{13}-l_{11}] 
\right)\times\\
&&\times
\left({[2l'][l_{33}-l_{31}]\over [2l+1][l_{11}-l_{23}]} 
\right)^{1/2}(m)_{13}^{-11+31+31}\\[2mm]
&   &\left. -{1\over a_{14}}{1\over [2l+1][2l'+1]} 
\left([2l][l_{13}-l_{11}][2l'][l_{33}-l_{31}] 
\right)^{1/2}(m)_{14}^{-11+31+31}\right\},\\[2mm] 
\left|1,1,1,1;(m)\right>& = &{1\over a_{15}}(m)^{-11-11+31+31}. 
\end{eqnarray*}
$$\eqno(A.1)$$
\end{document}